\documentclass[12pt]{article}
\usepackage{amsmath}
\usepackage{amssymb}
\usepackage{amsfonts}
\usepackage{graphics}

\newtheorem{thm}{Theorem}[section]

\newtheorem{prop}[thm]{Proposition}

\newenvironment{remark}{\par\medskip\noindent{\bf Remark.\ }}{\par\smallskip}

\newcommand{\be}{\begin{equation}}
\newcommand{\ee}{\end{equation}}
\newcommand{\openbox}{\leavevmode
  \hbox to8pt{\hfil\vrule\vbox to6pt{\hrule width6pt\vfil\hrule}\vrule}}

\newcommand{\qed}{\hbox to5pt{ } \hfill \openbox\bigskip\medskip}

\newcommand{\cF}{\mbox{$\cal F$}}
\newcommand{\cG}{\mbox{$\cal G$}}

\newcommand{\ve}[1]{\mathbf{#1}}

\newcommand{\R}{\mathbb R}

\newcommand{\F}{\mathbb F}

\title{Upper bounds for the size of ordered $L$-intersecting  set systems}
\author{G\'abor Heged\"{u}s
\\{\normalsize  \'Obuda University}
\\{\normalsize B\'ecsi \'ut 96/B, Budapest, Hungary, H-1032}
\\{\normalsize hegedus.gabor@uni-obuda.hu}
}

\begin{document}
\maketitle

\begin{abstract}
A family $\mbox{$\cal F$}=\{F_1,\ldots,F_m\}$  of subsets
of $[n]$  is said to be ordered, if there exists an $1\leq r\leq m$ index such that
$n\in F_i$ for each $1\leq i\leq r$,
$n\notin F_i$ for each $i>r$ and
$|F_i|\leq |F_j|$ for each $1\leq i<j\leq m$.

Our main result is a new upper bound for the size of ordered $L$-intersecting set systems.
\end{abstract}

\section{Introduction} 

Let  $[n]$ stand for the set $\{1,2,
\ldots, n\}$. We denote the family of all subsets of $[n]$  by $2^{[n]}$. 
Let $\binom{[n]}m$ denote the family of all subsets of $[n]$
which have cardinality $m$.

For a field $\F$, let 
$\F[x_1, \ldots, x_n]=\F[\ve x]$ denote  the
ring of polynomials in the variables $x_1, \ldots, x_n$ over $\F$.
For a subset $F \subseteq [n]$ we write
$\ve x_F = \prod_{j \in F} x_j$. In particular, $\ve x_{\emptyset}= 1$.

Let $\cF=\{F_1,\ldots,F_m\}$ be a family of subsets
of $[n]$. Let $L=\{\ell_1,\ldots ,\ell_s\}$ be a set of $s$ non-negative   integers. We say that the $\cF$ family is an $L$-intersecting family, if $|F_i \cap F_j|\in L$ for each $1\leq i,j\leq m, ~i \neq j$.

Frankl and Wilson proved the following general upper bound for the size of $L$-intersecting families.
\begin{thm} \label{FW}
Let $L=\{\ell_1,\ldots ,\ell_s\}$ be a set of $s$ non-negative integers. 
Let $\cF=\{F_1,\ldots,F_m\}$ be an $L$-intersecting  family of subsets
of $[n]$. Then 
$$
m\leq \sum_{i=0}^s {n \choose i}.
$$
\end{thm}

Snevily conjectured the following statement in his doctoral dissertation (see \cite{S3}).  Finally he proved this result in \cite{S2}.

\begin{thm} \label{Snevily}
 Let $L=\{\ell_1,\ldots ,\ell_s\}$ be a set of $s$ positive integers.  Let $\cF=\{F_1,\ldots,F_m\}$ be an $L$-intersecting  family of subsets
of $[n]$. Then   
$$
m\leq \sum_{i=0}^s {n-1 \choose i}.
$$
\end{thm}

We say that a set system $\cF=\{F_1,\ldots,F_m\}$ is {\em ordered}, if there exists an $1\leq r\leq m$ such that
\begin{enumerate}
\item 
$n\in F_i$ for each $1\leq i\leq r$;
\item
$n\notin F_i$ for each $i>r$ and
\item
$|F_i|\leq |F_j|$ for each $1\leq i<j\leq m$.
\end{enumerate}

Our main result is a new upper bound for the size of {\em ordered} $L$-intersecting set systems. Our proof is a version of the proof of Theorem \ref{Snevily}. 

\begin{thm} \label{main}
Let $L=\{\ell_1,\ldots ,\ell_s\}$ be a set of $s$ non-negative integer.  Let $\cF=\{F_1,\ldots,F_m\}$ be an ordered $L$-intersecting  family of subsets
of $[n]$. Then   
$$
m\leq \sum_{i=0}^s {n-1 \choose i}.
$$
\end{thm}
\begin{remark}
It is easy to verify that Theorem \ref{main} is sharp. Namely consider the set system $\cF:=\{A\subseteq [n]:~ n\notin A, |A|\leq s\}$. Then $\cF$ is an ordered $L$-intersecting  family of subsets, where $L=\{0,\ldots ,s-1\}$ and $|\cF|=\sum_{i=0}^s {n-1 \choose i}$.

We give an other example, which shows that Theorem \ref{main} is sharp: let 
$$
\cG:=\{G\subseteq [n]:~ n\in G, |G|\leq s\}\cup \{T\subseteq [n]:~ n\notin T, |T|= s\}.
$$ 
Now $\cG$ is an ordered $L$-intersecting  family of subsets, where $L=\{0,\ldots ,s-1\}$ and clearly $|\cG|=\sum_{i=0}^s {n-1 \choose i}$.
\end{remark}

\section{Proof}

The proof of  our main result is based on  the linear algebra bound method and the Triangular Criterion (see \cite{BF} Proposition 2.5). We recall here for the reader's convenience this principle.

\begin{prop} \label{tri} (Triangular Criterion)
Let $\F$ denote an arbitrary field. For $i=1,\ldots m$ let $f_i:\Omega \to \F$ be functions and $\ve v_i\in \Omega$ elements such that $f_i(\ve v_j)\neq 0$ if $i=j$ and $f_i(\ve v_j)=0$ if $j<i$. Then $f_1,\ldots ,f_m$ are linearly independent functions of the vector space $\F^{\Omega}$. 
\end{prop}

A polynomial is said to be {\em multilinear}, if it has degree at most $1$ in each variable. Let $f$ be a polynomial in $\F[x_1, \ldots, x_n]$ of degree at most $s$. Then there exists a unique  multilinear polynomial $\overline{f}$ of degree at most $s$ such that 
$$
f(\ve v)=\overline{f}(\ve v)
$$
for each $\ve v\in \{0,1\}^n$. This $\overline{f}$ polynomial is the  {\em multilinearization} of the polynomial $f$.

{\bf Proof of Theorem \ref{main}:}\\

Let $i\in [m]$ be a fixed index. Let $\ve v_i\in \{0,1\}^n$ denote the characteristic vector of the set $F_i$.

We denote by $\langle \ve x, \ve y\rangle:=\sum_{i=1}^n x_iy_i$ the usual scalar product of the vectors $\ve x$ and $\ve y$. Clearly $\langle \ve v_i, \ve v_j\rangle=|F_i\cap F_j|$.

For each $1\leq i\leq m$ let us define the real polynomials $P_i(x_1, \ldots ,x_n)$ as follows           
\begin{equation}  \label{Ppol} 
P_i(\ve x):= \prod_{k:\ell_k<|F_i|} \Big( \langle \ve x, \ve v_i\rangle-\ell_k \Big)\in \R[\ve x].
\end{equation}

It follows from the condition $|F_1|\leq \ldots \leq |F_m|$ that
$P_i(\ve v_i)\neq 0$ for each $1\leq i\leq m$ and $P_i(\ve v_j)=0$ for each $j<i$. 

Let $Q_i(\ve x)$ denote zhe multilinearization of $P_i(\ve x)$. Then $Q_i(\ve v_i)\neq 0$ for each $1\leq i\leq m$ and $Q_i(\ve v_j)=0$ for each $j<i$, hence Proposition \ref{tri} implies that the polynomials  $\{Q_i:~ 1\leq i\leq m\}$ are  linearly independent. Clearly  $\mbox{deg}(Q_i)\leq s$ for each $1\leq i\leq m$, because  $\mbox{deg}(P_i)\leq s=|L|$ for each $1\leq i\leq m$. 

Now we introduce $N:=\sum_{i=0}^{s-1} {n-1 \choose i}$ new  polynomials. Let $\cG:=\{T_1,\ldots ,T_N\}:=\cup_{i=0}^{s-1} {[n-1]\choose i}$  be the family of subsets  of $[n]$ with size at  most $s-1$, which doesn't contain $n$.

We can assume that $|T_1|\leq \ldots \leq |T_N|$. For each $1\leq i\leq N$ define the polynomial $g_i(\ve x):=(x_n-1)\cdot \prod_{j\in T_i} x_j$, where we understand $g_1(\ve x):=x_n-1$.  Let $\ve w_i\in \{0,1\}^n$ denote the characteristic vector of the set $T_i$  for each $1\leq i\leq N$.

It is easy to verify that $g_i(\ve w_i)\neq 0$ for each $1\leq i\leq m$ and $g_i(\ve w_j)=0$ for each $j<i$. Consequently Proposition \ref{tri} implies that the set of polynomials  $\{g_i:~ 1\leq i\leq N\}$ is linearly independent.

Now we prove that  the set of polynomials  $\{Q_i:~ r+1\leq i\leq m\}\cup\{g_i:~ 1\leq i\leq N\}$ is linearly independent.

Indirectly, if $\{Q_i:~ r+1\leq i\leq m\}\cup\{g_i:~ 1\leq i\leq N\}$ is not linearly independent, then there exists a non-trivial linear combination
\begin{equation}  \label{lincomb} 
\sum_{i=r+1}^m \alpha_i Q_i + \sum_{j=1}^N  \beta_j g_j=0. 
\end{equation}
Since we proved that the  polynomials $\{Q_i:~ r+1\leq i\leq m\}$ and  $\{g_j:~ 1\leq j\leq N\}$ are linearly independent, hence there exists  $r+1\leq i\leq m$ and $1\leq j\leq m$ indices such that $\alpha_i\neq 0$ and $\beta_j\neq 0$.

Let $j_0$ denote the minimal index such that  $\beta_{j_0}\neq 0$. Then clearly $\beta_j=0$ for each $1\leq j<j_0$. 

Now we determine the coefficient of the monomial $\ve x_{T_{j_0}}\cdot x_n$ in the  linear combination (\ref{lincomb}). 

It follows from the definition of the polynomial $Q_i$ that $Q_i\in \R[x_1, \ldots ,x_{n-1}]$ for each $r+1\leq i\leq m$, hence if we expand the polynomial $\sum_{i=r+1}^m \alpha_i Q_i$ as a sum of  monomials, then the monomial $\ve x_{T_{j_0}}\cdot x_n$ doesn't appear in this expansion. 

On the other hand we can expand also the sum $\sum_{j=1}^N  \beta_j g_j$ as a sum of  monomials. Since $\beta_j=0$ for each $1\leq j<j_0$, hence the only term in the sum $\sum_{j=1}^N  \beta_j g_j$, in which the monomial  $\ve x_{T_{j_0}}\cdot x_n$ appears, is  $\beta_{j_0} g_{j_0}$. Consequently the coefficient of the monomial $\ve x_{T_{j_0}}\cdot x_n$ in the  sum $\sum_{i=r+1}^m \alpha_i Q_i + \sum_{j=1}^N  \beta_j g_j$  is $\beta_{j_0}$, but this sum must be the zero  polynomial by  (\ref{lincomb}), hence $\beta_{j_0}=0$, a contradiction.

{\bf Claim.} The set of polynomials  $\{Q_i:~ 1\leq i\leq m\}\cup\{g_i:~ 1\leq i\leq N\}$ is linearly independent.

 {\bf Proof.} Consider the  linear combination 
\begin{equation}  \label{lincomb2} 
\sum_{i=1}^m \alpha_i Q_i + \sum_{j=1}^N  \beta_j g_j=0.
\end{equation}
Assume the contrary, that the  polynomials  $\{Q_i:~ 1\leq i\leq m\}\cup\{g_i:~ 1\leq i\leq m\}$ are not linearly independent. Since the polynomials  $\{Q_i:~ r+1\leq i\leq m\}\cup\{g_i:~ 1\leq i\leq m\}$ are linearly independent, hence we can suppose that $\alpha_{i_0}\neq 0$ for some $1\leq i_0\leq r$ and $\alpha_i=0$ for each $1\leq i< i_0$. But it follows from $n\in F_{i_0}$ that $g_j(\ve v_{i_0})=0$ for each $j$. Then by substituting $\ve v_{i_0}$ into  the equation  (\ref{lincomb2}) we get $\alpha_{i_0} Q_{i_0}(\ve v_{i_0})=0$, namely $Q_i(\ve v_{i_0})=0$ for each $i_0<i$. It follows from $Q_{i_0}(\ve v_{i_0})\ne 0$ that $\alpha_{i_0}=0$, a contradiction.

Let $V$ denote the vector space of multilinear polynomials in $n$ variables of degree at most $s$. Clearly $\dim V=\sum_{i=0}^s {n \choose i}$. 

We have found $m+N$  linearly independent  polynomials in $V$. Consequently we get
$$
|\cF|=m\leq \dim V - N=\sum_{i=0}^s {n \choose i} -\sum_{i=0}^{s-1} {n-1 \choose i}=\sum_{i=0}^s {n-1 \choose i}.
$$
\qed


\begin{thebibliography}{MM}

\bibitem{BF} L. Babai and P. Frankl, {\em Linear algebra methods in
combinatorics,} manuscript, September 1992.

\bibitem{FW} P. Frankl and  R.M. Wilson, Intersection theorems with geometric consequences. {\em Combinatorica} {\bf 1},  (1981) 357-368.

\bibitem{J} S. Jukna,  Extremal combinatorics: with applications in computer science. Springer Science and Business Media (2011).
 
\bibitem{S1} H. S. Snevily, On generalizations of the deBruijn-Erd\H{o}s theorem, {\em J. Comb. Theory
Ser. A} {\bf 68} (1994), 232-238.

\bibitem{S2}  H. S. Snevily,  A sharp bound for the number of sets that pairwise intersect at $k$ positive values. {\em Combinatorica} {\bf 23}(3), (2003) 527-533.


\bibitem{S3} H. S. Snevily,  Combinatorics of finite sets,  Doctoral dissertation, University of Illinois at Urbana-Champaign, 1991. 


\end{thebibliography}
\end{document}